\documentclass[12pt,reqno]{amsart}
\usepackage[a4paper]{geometry}
\usepackage{amsthm,hyperref,lmodern,mathrsfs,amssymb,amsmath,amsfonts}
\usepackage{setspace}
\usepackage[latin1]{inputenc}
\usepackage{versions}

\title[Concentration of invariant measures]{Measure concentration through non-Lipschitz observables and functional inequalities}

\author{Arnaud~Guillin} %
\address[A.~Guillin]{UMR CNRS 6620, Universit\'e Blaise Pascal, Clermont-Ferrand II and Institut Universitaire de France (IUF), France} %
\email{\url{mailto:guillin(at)math.univ-bpclermont.fr}} %
\urladdr{\url{http://math.univ-bpclermont.fr/~guillin/}}

\author{Ald\'eric~Joulin} %
\address[A.~Joulin, corresponding author]{UMR CNRS 5219, Institut de Math\'ematiques de Toulouse, Universit\'e de Toulouse, France}%
\email{\url{mailto:ajoulin(at)insa-toulouse.fr}}%
\urladdr{\url{http://www-gmm.insa-toulouse.fr/~ajoulin/}}



\numberwithin{equation}{section}

\newenvironment{Proof}{\removelastskip\par\medskip
\noindent{\em Proof.}
\rm}{\nolinebreak\hfill\rule{2mm}{2mm}\medbreak\par}

\newtheorem{theo}{Theorem}[section]
\newtheorem{lemme}[theo]{Lemma}

\newtheorem{defi}[theo]{Definition}
\newtheorem{remark}[theo]{Remark}

\def\Z{{\mathbb{Z}\ \!\!}}
\def\N{{\mathbb{N}\ \!\!}}
\def\R{{\mathbb{R}\ \!\!}}
\def\L{{\mathcal{L}\ \!\!}}
\def\D{{\mathcal{D}\ \!\!}}

\def\F{{\mathscr{F}\ \!\!}}
\def\B{{\mathscr{B}\ \!\!}}

\def\P{{\mathbb{P}\ \!\!}}

\def\cE{{\mathcal{E}\ \!\!}}
\def\X{{\mathcal{X}\ \!\!}}

\def\Var{{\mathrm{{\rm Var}}}}

\def\Ent{{\mathrm{{\rm Ent}}}}

\def\Cov{{\mathrm{{\rm Cov}}}}

\def\Lip{{\mathrm{{\rm Lip}}}}
\def\and{{\mathrm{{\rm and}}}}

\begin{document}

\begin{abstract}
Non-Gaussian concentration estimates are obtained for invariant probability measures of reversible Markov processes. We show that the functional inequalities approach combined with a suitable Lyapunov condition allows us to circumvent the classical Lipschitz assumption of the observables. Our method is general and covers diffusions as well as pure-jump Markov processes on unbounded spaces.
\end{abstract}

\maketitle

\section{Introduction}
\label{sect:intro} \textrm{} \indent
In the last few decades, the concentration of measure phenomenon has attracted a lot of attention. Given a metric probability space $(\X, d, \mu)$ and a sufficiently large class of functions defined on this space (we call them observables), the concentration of measure occurs when, observed through these functions, the space seems to be actually smaller than it is. In other words, there exists a non-decreasing continuous function $\alpha :[0,\infty) \to [0,\infty)$, null at the origin and tending to infinity at infinity, such that for a given (symmetric) class $\mathcal{C}$ of observables $f: \X \to \R$,
\begin{eqnarray*}
\mu \left(  \{ x\in \X : f(x) -\int_\X f \, d\mu  > r  \} \right) & \leq & e^{- \alpha (r)}, \quad r \geq 0.
\end{eqnarray*}
The concentration is said to be Gaussian when $\alpha$ is quadratic-like. In connection with isoperimetry theory, the class $\mathcal{C}$ is usually taken to be the space of Lipschitz functions on $(\X, d, \mu)$, say $\Lip (\X )$. A good review on the subject is the monograph \cite{ledoux} where the interested reader will find a clear introduction to the topic. One may mention also the recent progress in the area through mass transportation techniques, see the recent survey \cite{gozlan_leonard_survey}. \vspace{0.1cm}

In this paper, we emphasize a dynamical point of view of concentration of measure. Given the invariant measure $\mu$ of an ergodic continuous-time Markov process $(X_t)_{t\geq 0}$ with \textit{carr\'e du champ} operator $\Gamma$ (see below for the definition), we provide concentration properties of $\mu$ through observables which depend on the dynamics. As it is sometimes the case in previous studies, see for instance \cite{ledoux2} for a state of the art, our starting point is to assume that the pair $(\mu, \Gamma)$ satisfies a convenient functional inequality such as Poincar\'e or the entropic inequality. Such inequalities, which are verified by a wide variety of examples, are closely related to the long-time behaviour of the process. In particular, this approach allows to unify both continuous and discrete space settings even if, in essence, these two objects are rather different from each other. In the diffusion framework, the main ingredient to establish Gaussian concentration through functional inequalities is based on the chain rule derivation formula satisfied by the operator $\Gamma$. In the context of Markov jump processes, although this property is not verified - here $\Gamma$ becomes a finite difference operator - this difficulty can be circumvented by using reversibility and Gaussian concentration properties are still available in this discrete situation. In both cases the carr\'e du champ refers to a natural distance related to the dynamics and, within this notion of distance, the Lipschitz observables under which concentration estimates are obtained are the ones with a bounded carr\'e du champ, that is to say the space $\Lip_\Gamma (\X)$ of functions $f$ such that $\Gamma (f,f)$ is $\mu$-essentially bounded. Then a natural question arises: which type of measure concentration can we obtain beyond the space $\Lip_\Gamma (\X)$ ? In particular in discrete space settings, such a study basically makes sense on an unbounded state space $\X$. Using the notion of Ricci curvature for Markov chains (a definition also available in continuous-time), a first result of this kind was given by Ollivier in \cite{ollivier}, in which he obtains concentration bounds involving a mixed Gaussian-exponential regime, i.e. $\alpha (r)$ is quadratic/linear for small/large $r$. In our language, he requires that the carr\'e du champ $\Gamma(f,f)$ belongs to the space $\Lip (\X)$. Despite this interesting and new result, which is sufficiently robust to be extended to additive functionals, see e.g. \cite{joulin2} and \cite{jou_oll}, it seems to the authors that there is no satisfactory treatment yet to this question and we hope to give (the beginning of) an answer to this problem with the present article. \vspace{0.1cm}

Our idea is to use a Lyapunov condition on the observables. This kind of criteria have been successfully used for proving various types of functional inequalities, cf. \cite{bbcg, cgw, cgw2, glgy} and for concentration estimates of additive functionals, see for instance \cite{cg, ggw, cha_jou}. Namely we will consider the class $\L _V (a,b)$ of observables $f$ such that
$$
\Gamma (f,f) \leq -a \, \frac{\L V}{V} + b,
$$
where $a,b$ are two positive constants and $V$ is a convenient test function. When $a$ vanishes the class $\L _V (a,b)$ reduces to the space $\Lip_\Gamma (\X)$ and the classical concentration results apply \cite{ledoux2}. In particular, the behaviour of the carr\'e du champ $\Gamma(f,f)$ depends now on the growth of the term $-\L V /V$, which has no reason to be bounded. Of course there is a price to pay for such an improvement and it resides in the concentration property of the measure $\mu$, which is no longer Gaussian but only of Gaussian-exponential type under this class of observables. \\
To fix the ideas, let us consider a simple example on the set $\N := \{ 0,1,2, \ldots \}$ endowed with the classical distance between integers. Denote $\mu_p$ the geometric distribution on $\N$ of parameter $p \in (0,1)$, i.e. $\mu_p (\{ x \} ) := (1-p) p^x$, $x\in \N$. Then there exist infinitely many dynamics on $\N$ admitting $\mu_p$ as their reversible invariant measure, and among them let us consider the Markov process with carr\'e du champ
$$
\Gamma ^{(n)} (f,f) (x) := \frac{1}{2} \, \left \{ p(x+1)^n \, (f(x+1) - f(x)) ^2 + x^n 1_{\{x \neq 0\}} \, (f(x-1) - f(x)) ^2\right\},
$$
where $x\in \N$ and $n\in \N$ is some fixed parameter. On the one hand we can prove that for every $n\in \N$ the dynamics $(\mu_p,\Gamma ^{(n)})$ satisfies a Poincar\'e inequality and adapting then to the discrete case the method introduced by Aida and Stroock in \cite{aida_stroock} entails a mixed Gaussian-exponential regime for the measure $\mu_p$, under observables $f\in \Lip_{\Gamma ^{(n)}}(\N)$. Actually, these observables correspond to:

- functions belonging to the space $\Lip (\N)$ in the case $n=0$. Such an estimate is sharp, cf. for instance \cite{houdre, joulin0};

- functions of order $\sqrt{x}$ in the case $n=1$. According to the previous line, the measure $\mu_p$ should concentrate at least like a Gaussian through these observables;

- functions of order $\log (x)$ in the case $n=2$. Under these observables, $\alpha (r) $ should behave like ${e^{-e^r}}$ for large $r$;

- bounded functions in the case $n\geq 3$. \\
On the other hand Ollivier's method cleverly applies, but only in the case $n=1$. Hence one observes that these two methods does not lead to sharp results as soon as $n\geq 2$. Actually, thanks to the Lyapunov approach, the Gaussian-exponential concentration results we give in this paper apply for observables $f$ with carr\'e du champ $\Gamma ^{(n)} (f,f)(x)$ of order $x^n$, i.e. of the expected order of magnitude. \vspace{0.1cm}

The paper is organized as follows. In section~\ref{sect:nota_ergodic}, we recall some basic material on Markov processes and functional inequalities. Two types of processes are considered in our study: diffusions and pure-jump Markov processes. Next we state in section~\ref{sect:result} our main results of the paper, theorems~\ref{theo:result} and \ref{theo:poincare}, in which some mixed Gaussian-exponential concentration properties of $\mu$ are obtained through observables satisfying the Lyapunov condition defined above and under a convenient functional inequality assumption satisfied by the dynamics $(\mu,\Gamma)$. As a result, such new concentration inequalities extend the classical estimates obtained when the observables belong to the space $\Lip_\Gamma (\X)$, corresponding to the case where $a$ vanishes. Finally, section~\ref{sect:ex} is devoted to numerous examples in continuous and discrete settings. First we investigate the diffusion case of Kolmogorov processes whose invariant measure has density proportional to $e^{-U}$, where $U$ is some nice potential on $\R^d$, and for which we derive a concentration result when $U$ is considered as an observable, in the spirit of the recent progress made by Bobkov and Madiman in \cite{bob_madiman}. Dealing with jump processes, the case of birth-death processes is addressed in detail. When we apply our concentration estimates to observables in $\Lip (\N)$ equipped with the classical distance between integers, we are able to consider processes whose generator is unbounded, like the basic example investigated above. Finally, we focus our attention on an unbounded interacting particle system, namely the Glauber dynamics associated to a Gibbs measure defined with respect to a Poisson reference measure.

\section{Preliminaries}
\label{sect:nota_ergodic}

\subsection{Functional inequalities}
Throughout the paper, $(\X,d)$ is a Polish space endowed with the corresponding Borel $\sigma$-field $\B$ and $\Lip (\X )$
is the space of Lipschitz functions on $\X$ with finite Lipschitz seminorm with respect to $d$, i.e.
\begin{eqnarray*}
\Vert f\Vert _{\Lip } & := & \sup _{x\neq y } \frac{\vert f(x)-f(y) \vert }{d(x,y)} <\infty .
\end{eqnarray*}
On a filtered probability space $(\Omega , \F , (\F_t)_{t\geq 0}, \P )$, let $\left \{ (X_t)_{t\geq 0}, (\P _x )_{x\in \X } \right \}$ be an
$\X$-valued \textit{c\`adl\`ag} ergodic Markov process with reversible invariant measure (or stationary distribution) $\mu$ and symmetric semigroup $(P_t)_{t\geq 0}$ on $L^2(\mu)$. In the sequel we denote $L^p(\mu) := L^p(\X , \B , \mu)$ for $p\in [1,\infty]$. Denote $\L$ the self-adjoint generator acting on its dense domain $\D _2 (\L)$ consisting of functions $f \in L^2 (\mu )$ such that $t^{-1}(P_t f - f)$ admits a limit in $L^2(\mu)$ as $t\to 0$. One of the main protagonists of the present paper is the \textit{carr\'e du champ} $\Gamma$, which is a bilinear symmetric operator
on $\D _2 (\L) \times \D _2 (\L)$ defined by
\begin{eqnarray*}
\Gamma (f,g) & := & \frac{1}{2} \, \left( \L (fg) - f \L g  - g  \L f \right) .
\end{eqnarray*}
As mentioned in the introduction, there is a natural pseudo-distance associated to the
operator $\Gamma$ which can be defined as
$$
d_\Gamma (x,y) : = \sup \left \{ \vert f(x)-f(y) \vert : \Vert \Gamma(f,f) \Vert _{L^\infty (\mu)} \leq 1 \right \} , \quad x,y\in \X  .
$$
Although this distance can be infinite, it is well-defined in the situations of interest and carries
a lot of information about the structure of the underlying process. In the sequel, we denote $\Lip_\Gamma (\X) $ the space of Lipschitz functions with respect to $d_\Gamma$. \\
The associated Dirichlet form acts on $\D (\cE _\mu) \times \D (\cE _\mu)$ as
\begin{eqnarray*}
\cE _\mu (f,g) & = & \int _\X \Gamma (f,g) \, d\mu ,
\end{eqnarray*}
where $\D (\cE _\mu)$ is the subspace of all functions $f \in \D _2 (\L)$ such that $\cE _\mu (f,f)$ is well-defined. In particular,
the Donsker-Varadhan information of any probability measure $\nu$ on $\X$ with respect to the invariant measure $\mu$ is defined as
$$
I(\nu | \mu ) :=
\left \{
\begin{array}{ll}
\cE _\mu \left( \sqrt{f},\sqrt{f} \right) & \mbox{if} \quad d\nu = f d\mu , \quad \sqrt{f} \in \D (\cE _\mu) ; \\
\infty & \mbox{otherwise}.
\end{array}
\right.
$$
A key point in our analysis is that the functional $\nu \mapsto I(\nu |\mu)$ is nothing but the rate function governing the Large Deviation
Principle in large time of the empirical measure of $(X_t)_{t\geq 0}$. However in the non-reversible case, it is
given by a contraction form of the Donsker-Varadhan entropy which is different from the Donsker-Varadhan information, so that
our study will not extend to the non-symmetric case, unfortunately. \vspace{0.1cm}

Now let us introduce the functional inequalities we will focus on in the paper. Given an integrable function $f\in L^1(\mu)$,
we denote $\mu (f) := \int _\X f d\mu$. Let $I$ be an open interval of $\R$ and for a convex function $\phi: I \to \R$ we define
the $\phi$-entropy of a function $f : \X \to I$ with $\phi(f) \in L^1(\mu)$ as
$$
\Ent_\mu ^{\phi} (f) := \mu \left( \phi(f) \right) - \phi \left( \mu (f) \right).
$$
The dynamics $(\mu, \Gamma)$ satisfies a $\phi$-entropy inequality with
constant $C_\phi > 0$ if for any $I$-valued function $f\in \D (\cE _\mu)$ such that $\phi'(f) \in \D (\cE _\mu)$,
$$
C_\phi \, \Ent_\mu ^{\phi} (f) \leq \frac{1}{2} \, \cE _\mu (f, \phi '(f)).
$$
See for instance \cite{chafai} for a careful investigation of the properties of $\phi$-entropies. The latter inequality is satisfied if and only if the following entropy dissipation of the semigroup holds: for any $I$-valued function
$f$ such that $\phi(f) \in L^1(\mu)$,
$$
\Ent_\mu ^\phi (P_tf) \leq e^{-2 C_\phi t } \, \Ent_\mu ^\phi (f), \quad t\geq 0.
$$
In this paper we will consider three cases: \vspace{0.1cm}

$(i)$ the Poincar\'e inequality: $\phi(u) = u^2$ with $I=\R$ and the $\phi$-entropy inequality rewrites as
\begin{equation}
\label{eq:poincare}
\lambda \, \Var _\mu (f) \leq \cE _\mu \left( f,f\right) ,
\end{equation}
where the variance of $f$ under $\mu$ is given by
$$
\Var _\mu (f) := \mu (f^2) - \mu (f) ^2.
$$
The optimal constant $\lambda _1$ (say) is nothing but the spectral gap in $L^2(\mu)$ of
the operator $-\L$, i.e. its smallest non-zero eigenvalue. Estimating $\lambda _1$ allows us to obtain the optimal
rate of convergence of the semigroup in $L^2(\mu)$. \vspace{0.1cm}

$(ii)$ the entropic inequality: $\phi(u) = u \, \log u$ with $I=(0,\infty)$ and the $\phi$-entropy inequality is given by
\begin{equation}
\label{eq:mlsi}
\rho \, \Ent _\mu (f) \leq \cE _\mu \left( f,\log f\right) ,
\end{equation}
where the entropy under $\mu$ of the smooth positive function $f$ is defined by
$$
\Ent _\mu (f) := \mu (f \, \log f ) - \mu (f) \, \log \mu (f) .
$$
We have skipped in the inequality the constant $1/2$ for convenience in future computations. Once again, the best constant
$\rho _0$ in (\ref{eq:mlsi}) gives the optimal exponential decay of the entropy along the semigroup. \vspace{0.1cm}

$(iii)$ the Beckner-type inequality: $\phi(u) = u^p $ with $p\in (1,2]$ and $I=(0,\infty)$. We have in this case
\begin{equation}
\label{eq:tetali}
\alpha_p \, \left( \mu (f^p) - \mu(f) ^p \right) \leq \frac{p}{2} \, \cE _\mu \left( f,f^{p-1}\right) .
\end{equation}
Estimating $\alpha_p$ gives the optimal rate of convergence of the semigroup in $L^p(\mu)$. \vspace{0.1cm}

The entropic and Beckner-type inequalities are stronger than the Poincar\'e inequality (apply these
inequalities to the function $1+\varepsilon f$ and take the limit as $\varepsilon \to 0$). Moreover it reduces to the
Poincar\'e inequality (\ref{eq:poincare}) if $p=2$, whereas dividing both sides by $p-1$
and taking the limit as $p\to 1$ we obtain the entropic inequality (\ref{eq:mlsi}). \vspace{0.1cm}

In this paper we will mainly consider two general classes of reversible Markov processes: diffusions and pure jump Markov processes to which we turn now.

\subsection{Diffusion processes}
\label{sous-sect:diff}
A diffusion process on the Euclidean space $\X = \R^d$ corresponds to a path continuous Markov process on $\R^d$ whose generator $\L$ is a
second order differential operator: for any sufficiently smooth function $f:\R^d \to \R$,
$$
\L f (x) = \sum_{i,j=1}^d a_{i,j} (x) \, \frac{\partial ^2 f}{\partial x_i \partial x_j} (x) + \sum_{i=1}^d b_i (x) \,
\frac{\partial  f}{\partial x_i } (x), \quad x\in \R^d .
$$
Here $a := \sigma \sigma ^*$ is a measurable and locally bounded function from $\R^d$ to the space of $d\times d$ symmetric positive definite matrices with smooth entries, $\sigma^*$ being the transpose of the matrix $\sigma$, and the measurable drift $b:\R^d \to \R^d$ is also assumed to be smooth. In this case the carr\'e  du champ is given by
\begin{eqnarray*}
\Gamma (f,g) & = & \sum_{i,j=1}^d a_{i,j} \, \frac{\partial f}{\partial x_i} \, \frac{\partial f}{ \partial x_j} \\
& = & < \sigma^* \nabla f , \sigma^* \nabla g > ,
\end{eqnarray*}
where $<\cdot , \cdot >$ stands for the Euclidean scalar product in $\R^d$ and $\nabla$ is the usual gradient operator. In particular when $\sigma$ is the identity matrix, the spaces $\Lip (\R^d)$ and $\Lip_\Gamma (\R^d)$ might be identified. \\
In contrast to the jump case introduced below, $\Gamma$ is a differentiation, i.e. for all
any smooth enough functions $(f_k)_{1\leq k \leq n}, f : \R^d \to \R$ and any $C^1$ function $\phi :\R^n \to \R$,
\begin{eqnarray}
\label{eq:chain_rule}
\Gamma (\phi(f_1, \ldots , f_n),f) = \sum_{i=1}^n \frac{\partial \phi}{\partial x_i} (f_1, \ldots , f_n) \, \Gamma (f_i, f).
\end{eqnarray}
Due to this chain rule derivation formula, the entropic inequality (\ref{eq:mlsi}) rewrites in the diffusion case as the
famous log-Sobolev inequality
\begin{equation}
\label{eq:lsi}
\rho \, \Ent _\mu (f^2) \leq 4 \, \cE_\mu (f,f) ,
\end{equation}
which is the original inequality (up to the extra factor 4) derived by Gross in \cite{gross} to study hypercontractivity of the underlying
semigroup. When we will consider diffusion processes in the sequel, we will use the terminology ``log-Sobolev inequality" instead of ``entropic inequality". \\
On the other hand, letting $p=2/q$ for $q\in [1,2)$ and $f=g^q$, the Beckner-type inequality (\ref{eq:tetali})
rewrites as the so-called standard Beckner inequality:
\begin{equation}
\label{eq:beckner}
\alpha _{2/q} \, \left( \mu (g^2) - \mu (g^{q}) ^{2/q} \right) \leq (2-q) \, \cE _\mu (g,g).
\end{equation}
Such an inequality was introduced by Beckner in \cite{beckner} for the Gaussian measure. In particular, the limiting case $q\to 2$ recovers the classical log-Sobolev inequality. Note however that the inequality (\ref{eq:beckner}) is weaker than the log-Sobolev inequality, cf. \cite{latala}.

\subsection{Markov jump processes}
\label{sous-sect:jump}
Dealing with a pure-jump Markov process, the generator $\L$ is given for any function $f \in \D _2 (\L)$ by
$$
\L f(x) = \int _\X \left( f(y)-f(x) \right) Q_x (dy), \quad x\in \X ,
$$
where the transition kernel $x\mapsto Q_x$ is a measurable mapping from $\X$ to the set of Radon measures on $\X$ endowed with the
corresponding Borel $\sigma$-field. We assume that it satisfies the following
stability assumption:
\begin{equation}
\label{eq:stab}
\int _\X Q_x (dy) <\infty , \quad x\in \X ,
\end{equation}
which entails that the process is piecewise constant. Here, reversibility means that the following detailed balance condition is satisfied:
\begin{equation}
\label{eq:revers}
Q_x (dy) \, \mu (dx) = Q_y (dx) \, \mu (dy) .
\end{equation}
The carr\'e du champ operator $\Gamma$ admits an explicit expression given for any $f,g \in \D _2 (\L)$ by
$$
\Gamma (f,g) (x) = \frac{1}{2} \, \int _\X \left( f(y) - f(x) \right) \left( g(y) - g(x) \right)  \, Q_x (dy) ,
$$
and we have
$$
\Gamma (f,f) (x) = \frac{1}{2} \, \int _\X \left( f(y) - f(x) \right) ^2 \, Q_x (dy) .
$$
In particular, the spaces $\Lip (\X)$ and $\Lip _\Gamma (\X)$ have no reason to coincide since the kernel of the generator may be unbounded, i.e.
\begin{equation}
\label{eq:unbounded}
\sup _{x\in \X} \, \int _\X Q_x (dy) =  \infty .
\end{equation}
Finally the Dirichlet form rewrites for any $f,g\in \D (\cE _\mu)$ as
\begin{eqnarray*}
\cE _\mu (f,g) & = & \frac{1}{2} \,  \int _\X \int _\X \left( f(y) - f(x) \right) \left( g(y) - g(x) \right) \, Q_x (dy) \mu (dx) \\
& = & \, \underset{f(x)>f(y)}{ \int \int } \left( f(y) - f(x) \right) \left( g(y) - g(x) \right) \, Q_x (dy) \mu (dx),
\end{eqnarray*}
where in the last line the reversibility is used. In our jump framework, the entropic inequality (\ref{eq:mlsi}) corresponds to one of the so-called modified log-Sobolev inequalities
introduced by Bobkov and Ledoux in \cite{bob_ledoux}. However, due to the lack of chain rule for discrete gradients, this inequality
is different from the discrete version of the log-Sobolev inequality (\ref{eq:lsi}), and the same remark holds between the Beckner-type
inequality (\ref{eq:tetali}) and the standard Beckner inequality (\ref{eq:beckner}). We refer to \cite{diaco_saloff, bob_ledoux, bob_tetali}
for historical and tutorial references on these discrete functional inequalities, together with a hierarchy of the various modified
log-Sobolev inequalities.

\section{Main results}
\label{sect:result}
As announced, we obtain concentration properties of the invariant measure $\mu$ through observables which are not required to belong to the spaces $\Lip (\X)$ nor $\Lip _\Gamma (\X)$, but which satisfy a Lyapunov condition. In order to state this condition properly, let us introduce first the extended domain of the generator. Denote the probability measure $\P_\nu (\cdot ) := \int_\X \P_x (\cdot ) \, \nu (dx) $ where $\nu$ is an arbitrary initial probability distribution.
A continuous function $f$ is said to belong to the extended domain $\D_e (\L)$ of the
generator $\L$ if there exists some measurable function $g:\X \to \R$ such that for any $t\geq 0$, $\int_0 ^t \vert g(X_s) \vert \, ds < \infty$,
$\P_\mu$-a.s. and the process
$$
M_t ^f = f(X_t) - f(X_0) - \int_0 ^t g(X_s) \, ds , \quad t\geq 0,
$$
is a local $\P_\mu$-martingale. In this case we write $f\in \D_e (\L)$ and $\L f = g$. \vspace{0.1cm}

The first result on which our analysis is based is closely related to the theory of large deviations, see for instance \cite{wu}.
\begin{lemme}
\label{lemme:fisher}
The Donsker-Varadhan information $\nu \mapsto I(\nu |\mu)$ is the rate function governing the Large Deviation Principle of the
empirical measure of the reversible process $(X_t)_{t\geq 0}$. In other words,
\begin{equation}
\label{eq:LDP}
I(\nu |\mu) = \sup _{0 < V\in \D _e (\L)} \, \int _\X - \frac{\L V}{V} \, d \nu .
\end{equation}
\end{lemme}
\begin{remark}
\rm{Actually, we will only need the following inequality, available for any positive test function $V\in \D _e (\L)$:
$$
\int _\X - \frac{\L V}{V} \, d \nu \, \leq \, I(\nu |\mu) .
$$
For instance in the diffusion case, we can prove it as follows: assume that $\nu$ has density $f^2$ with respect to $\mu$ (trivial otherwise). Then by the chain rule formula satisfied by the carr\'e du champ,
\begin{eqnarray*}
\int _\X - \frac{\L V}{V} \, d \nu & = & \int _\X \Gamma \left( V, \frac{f^2}{V} \right) \, d \mu \\
& = & \int _\X \left( \frac{2f}{V} \, \Gamma (V,f)- \frac{f^2}{V^2} \, \Gamma (V,V) \right) \, d \mu \\
& \leq & \int _\X \left( \frac{2f}{V} \, \sqrt{\Gamma (V,V)} \, \sqrt{\Gamma (f,f)} - \frac{f^2}{V^2} \, \Gamma (V,V) \right) \, d \mu \\
& \leq & \int _\X \Gamma \left( f,f \right) \, d \mu \\
& = & I(\nu |\mu) ,
\end{eqnarray*}
where we used Cauchy-Schwarz' inequality to obtain the first inequality, and also the elementary bound $2xy \leq x^2 + y^2$, $x,y \in \R$, to get the final result.}
\end{remark}

Now we are able to state the Lyapunov condition we will focus on along this paper.
\begin{defi}
\label{hyp:jumps}
Let $a,b $ be two positive constants and let $V \in \D _e (\L)$ be a positive test function. A function $f\in \D _2 (\L)$
belongs to the class $\L _V(a,b)$ if the following inequality is satisfied $\mu-a.s.$:
\begin{equation}
\label{eq:lyap}
\Gamma (f,f)  \leq - a \, \frac{\L V}{V} + b.
\end{equation}
\end{defi}
\begin{remark}
\label{rem:poincare}
\rm{The Poincar\'e inequality can be seen as a minimal assumption in our study of concentration by means of the Lyapunov condition (\ref{eq:lyap}). Indeed, if there exists a function $f\in \D _2 (\L)$ such that $\Gamma (f,f)$ is lower bounded by a positive constant at infinity, and this the case in the main examples of interest (except in the Cauchy-like case appearing in section~\ref{sect:ex}), then the Poincar\'e inequality is satisfied, cf. \cite{cat_guil_zitt}. Moreover, integrating with respect to $\mu$ entails that $f\in \D (\cE _\mu)$ and by Poincar\'e inequality we have $\Var _\mu (f) \leq b /\lambda _1$. In other words the constant $b/\lambda_1$ can be interpreted in the sequel as a variance term of the observable $f$.}
\end{remark}
Before stating our first main result, let us provide a key lemma. In the remainder of this paper, we only give the proofs in the jump case since the diffusion framework requires no additional difficulties and is even simpler, according to the chain rule derivation formula (\ref{eq:chain_rule}) satisfied by the carr\'e du champ.
\begin{lemme}
\label{lemme:cle}
Let $f$ belong to the class $\L _V(a,b)$. Given $\lambda \in (0, 2/\sqrt{a})$, let $\mu_\lambda$ be the probability measure with density
$f_\lambda := e^{\lambda f} /Z_\lambda $ with respect to $\mu$, where $Z_\lambda$ is the appropriate normalization constant. We assume that $\sqrt{f_\lambda} \in \D (\cE _\mu)$. Then
$$
I(\mu_\lambda | \mu) \leq \frac{\lambda^2 b}{4-\lambda^2 a}, \quad 0<\lambda < \frac{2}{\sqrt{a}}.
$$
\end{lemme}
\begin{Proof}
Since $f\in \L _V (a,b)$, we have for any $\lambda \in (0, 2/\sqrt{a})$:
\begin{eqnarray*}
I(\mu _\lambda |\mu)
& = & \frac{1}{ Z_\lambda } \, \underset{f(x)>f(y)}{ \int \int }\left( e^{\lambda f(x) /2} -  e^{\lambda f(y) /2}\right) ^2 \, Q_x (dy) \mu (dx) \\
& = & \underset{f(x)>f(y)}{ \int \int } \left( 1- e^{-\lambda (f(x)-f(y))/2} \right) ^2 \, f_\lambda (x) \, Q_x (dy) \mu (dx) \\
& \leq & \frac{\lambda^2}{4} \, \int _\X \Gamma (f,f) \, d\mu _\lambda \\
& \leq & \frac{\lambda^2}{4} \, \int _\X \left( - a \, \frac{\L V}{V} + b \right) \, d\mu _\lambda \\
& \leq & \frac{\lambda ^2}{4} \, \left( a I(\mu _\lambda |\mu) + b\right) ,
\end{eqnarray*}
where in the last line we used lemma~\ref{lemme:fisher}. Finally rearranging the terms allows us to obtain the desired inequality.
\end{Proof}
We turn now to our first main and new result which exhibits a non-Gaussian concentration estimate through observables belonging to the class $\L _V(a,b)$. Due to the approach we will use, the numerical constants in the estimates below have no reason to be sharp.
\begin{theo}
\label{theo:result}
Assume that the pair $(\mu,\Gamma)$ satisfies the entropic inequality (\ref{eq:mlsi}) with optimal constant $\rho _0$. Let $f \in \L _V(a,b)$ and let
$$
r_{\max} := \frac{8b}{3\rho_0 \sqrt{a}}
$$
be the size of the Gaussian window. Then the invariant measure $\mu$ has the following concentration property: for any
$0\leq r \leq r_{\max} $, the deviation is of Gaussian-type:
\begin{equation}
\label{eq:conc2}
\mu \left( \left \{ x\in \X: f(x) - \mu (f) >r \right \} \right) \leq e^{ - \frac{3 \rho _0 r^2 }{16b} } ,
\end{equation}
and for any $r\geq r_{\max}$, the decay is exponential:
\begin{equation}
\label{eq:conc21}
\mu \left( \left \{ x\in \X: f(x) - \mu (f) >r \right \} \right) \leq e^{ -  \frac{r}{2\sqrt{a}}} ,
\end{equation}
\end{theo}
\begin{remark}
\rm{In the sequel, a concentration property such as (\ref{eq:conc2})-(\ref{eq:conc21}) will be called Gaussian-exponential concentration.}
\end{remark}
\begin{proof}
Denote $L_\lambda := \lambda ^{-1} \log Z_\lambda$, where $Z_\lambda := \int _\X e^{\lambda f} \, d\mu$, with $\lambda \in (0,1/\sqrt{a})$, and let $\mu _\lambda $ be the absolutely continuous probability measure with density $f_\lambda := e^{\lambda f}/ Z_\lambda$ with respect to $\mu$. Using a standard approximation procedure one may assume that the observable $f \in \L _V (a,b)$ is bounded so that $\sqrt{f_\lambda } \in \D (\cE _\mu)$. The following proof is a modification of the famous Herbst method popularized by Ledoux. Using the entropic inequality (\ref{eq:mlsi}),
\begin{eqnarray*}
\label{eq:detail}
\nonumber \frac{d}{d\lambda} L_\lambda & = & \frac{1}{\lambda ^2 Z_\lambda } \, \Ent _\mu (e^{\lambda f}) \\
\nonumber & \leq & \frac{1}{\rho _0 \lambda ^2 Z_\lambda} \, \cE _\mu \left( \lambda f,e^{\lambda f}\right) \\
\nonumber & = & \frac{1}{\rho _0 \lambda } \, \underset{f(x)>f(y)}{ \int \int }\left( f(x)-f(y) \right) \left( 1- e^{- \lambda (f(x)-f(y))} \right)
\, f_\lambda (x) \, Q_x (dy) \mu (dx) \\
\nonumber & \leq & \frac{1}{\rho _0 } \, \int _\X \int _\X  \Gamma (f,f) \, d\mu _\lambda \\
& \leq & \frac{1}{\rho _0} \, \int _\X \left( - a \frac{\L V}{V} + b\right) \, d\mu _\lambda \\
& \leq & \frac{1}{\rho _0} \, \left( a I(\mu_\lambda |\mu ) + b\right),
\end{eqnarray*}
where we used that $f\in \L _V(a,b)$ and then lemma~\ref{lemme:fisher} in the two last lines. Thus lemma~\ref{lemme:cle} entails the inequality
\begin{eqnarray*}
\frac{d}{d\lambda} L_\lambda & \leq & \frac{4b}{3\rho _0 } , \quad 0< \lambda < \frac{1}{\sqrt{a}},
\end{eqnarray*}
and therefore the following log-Laplace estimate is available for any $0< \lambda < \frac{1}{\sqrt{a}}$:
\begin{equation}
\label{eq:laplace}
\log \int _\X e^{\lambda f } \, d\mu \leq \lambda \mu (f) + \frac{4 b \lambda^2 }{3\rho _0} .
\end{equation}
Finally using Chebyshev's inequality and optimizing in $\lambda \in (0,1/\sqrt{a})$ yields the tail estimates (\ref{eq:conc2}) and (\ref{eq:conc21}). The proof of theorem~\ref{theo:result} is thus complete.
\end{proof}
\begin{remark}\rm{Two deviation regimes appear, Gaussian and exponential, with continuous transition from one to the other. In contrast to the classical Herbst method where the observables belong to $\Lip_\Gamma (\X)$, i.e. $a=0$ in the Lyapunov condition (\ref{eq:lyap}), our assumption allows us to go beyond this Lipschitz property. However the price to pay is to have a finite Gaussian window, i.e. $r_{\max} <\infty$.}
\end{remark}
\begin{remark}\rm{By the Central Limit Theorem, the order of magnitude is correct in the Gaussian regime. Since the entropic inequality entails a Poincar\'e inequality, we have $\rho _0 \leq \lambda _1$ and thus for any observable $f\in \L _V (a,b)$, we get $\Var_\mu (f) \leq b/\rho _0$. Therefore, if $\mu = \nu ^{\otimes d}$ is a product measure and $f(x) = \sum_{k=1}^d \phi(x_k)$, $x:=(x_1, \ldots, x_d) \in \X$, then we obtain the following inequality, which is sharp for large $d$:
$$
\mu \left( \left \{ x\in \X: f(x) - \mu (f) >r \sqrt{d} \right \} \right)
\leq e^{ - \frac{3 \rho _0 r^2 }{16\tilde{b}} } , \quad 0\leq r \leq \frac{8\tilde{b}\sqrt{d}}{3\rho_0 \sqrt{a}}.
$$
Here the important point is that the positive parameter $\tilde{b}$ depends on $\phi$ but is independent of $d$.}
\end{remark}
\begin{remark} \label{rem:regimes}\rm
The method is sufficiently robust to get, for large deviation level $r$, other regimes than exponential under particular observables. For example, assume that we consider $f\in \L _V(a,b)$ such that $\Gamma(f,f)\ll -a\L V/V+b$ but that there exists two functions $\phi,\psi : (0,\infty) \to (0,\infty)$ such that for all $\varepsilon>0$,
$$\Gamma(f,f)\le \phi(\varepsilon)\left(-a \, \frac{\L V}{V} +b\right)+\psi(\varepsilon).$$
Then plugging this estimate in the previous proofs of lemma \ref{lemme:cle} and theorem \ref{theo:result}, one has for all $\varepsilon >0$,
$$
I(\mu_\lambda|\mu) \le \frac{\lambda^2}{4}\left(a\phi(\varepsilon) I(\mu_\lambda|\mu) +b\phi(\varepsilon)+\psi(\varepsilon)\right).
$$
Optimizing in $\varepsilon >0$ enables to get for some function $\Phi : [0,\infty) \to [0,\infty ]$ and all $\lambda >0$,
$$
I(\mu_\lambda|\mu)\le \Phi(\lambda) ,
$$
leading then to super-exponential regime for large $r$. We will illustrate this on an example in section 4.
\end{remark}

Inspired by Otto-Villani's method, cf. \cite{otto_villani, sammer} where the links between log-Sobolev and transportation inequalities are studied on continuous and finite state spaces respectively, let us recover theorem~\ref{theo:result} by using a semigroup proof. Once again we focus our attention on the jump case. Let $h$ be a smooth density with respect to $\mu$. Given $t>0$, denote $\nu_t$ the probability measure with density $P_t h$ with respect to $\mu$. We assume that the Donsker-Varadhan information $I(\nu_t |\mu)$ is well-defined, i.e. $\sqrt{P_t h} \in \D (\cE _\mu)$. Using Cauchy-Schwarz's inequality and then reversibility,
\begin{eqnarray*}
\cE _\mu (P_t h , f) & = & \frac{1}{2} \, \int _\X \int_\X \left( P_t h(x) - P_t h(y)\right) \, \left( f(x)-f(y) \right) \, Q_x (dy) \mu (dx) \\
& \leq & \sqrt{I(\nu_t |\mu)} \, \sqrt{\frac{1}{2} \, \int _\X \int_\X \left( \sqrt{P_t h(x)}+ \sqrt{P_t h(y)} \right) ^2 \, \left( f(x)-f(y) \right) ^2 \, Q_x (dy) \mu (dx)} \\
& \leq & 2 \, \sqrt{I(\nu_t |\mu)} \, \sqrt{\int_\X \Gamma (f,f) \, d\nu_t } \\
& \leq & 2 \, \sqrt{I(\nu_t |\mu)} \, \sqrt{\int_\X \left( -a \, \frac{\L V}{V} + b \right) \, d\nu_t } \\
& \leq & 2 \, \sqrt{I(\nu_t |\mu)} \, \sqrt{a I(\nu_t |\mu)+ b} ,
\end{eqnarray*}
where in the two last lines we used that $f\in \L _V(a,b)$ and then lemma~\ref{lemme:fisher}. Using now the elementary inequality $2(a-b)^2 \leq (a^2 -b^2) \log (a/b)$ available for any $a,b>0$ and then the entropic inequality (\ref{eq:mlsi}), we get
\begin{eqnarray*}
\cE _\mu (P_t h , f) & \leq & \cE _\mu (P_t h , \log P_t h ) \, \left( \frac{\sqrt{a}}{2} + \sqrt{\frac{b}{\cE _\mu (P_t h , \log P_t h )}} \right) \\
& \leq & \cE _\mu (P_t h , \log P_t h ) \, \left( \frac{\sqrt{a}}{2} + \sqrt{\frac{b}{\rho_0 \Ent _\mu (P_t h )}} \right) .
\end{eqnarray*}
Integrating time between 0 and infinity entails the covariance inequality:
\begin{eqnarray*}
\Cov _\mu (f,h) & := & \mu (fh) - \mu (f) \, \mu (h) \\
& = & \int_0 ^\infty \cE _\mu (P_t h , f) \, dt \\
& \leq & \int_0 ^\infty \cE _\mu (P_t h , \log P_t h ) \, \left( \frac{\sqrt{a}}{2} + \sqrt{\frac{b}{\rho_0 \Ent _\mu (P_t h )}} \right) \, dt \\
& = & \frac{\sqrt{a}}{2} \, \Ent _\mu (h) + 2 \, \sqrt{\frac{b \, \Ent_\mu (h) }{\rho_0 }} ,
\end{eqnarray*}
which in turn yields the inequality
\begin{eqnarray*}
\alpha \left( \Cov _\mu (f,h)\right) & \leq & \Ent_\mu (h),
\end{eqnarray*}
where $\alpha$ is the function
$$
\alpha (r) = \frac{\rho_0 r^2}{4b + 2\rho_0 \sqrt{a}\, r}, \quad r>0.
$$
Finally using theorem~3.2 in \cite{gozlan_leonard_survey}, we obtain the following concentration estimate through the observable $f \in \L _V (a,b)$:
\begin{eqnarray*}
\mu \left( \left \{ x\in \X : f(x) -\mu (f)  > r \right \} \right) & \leq & e^{- \alpha (r)}, \quad r \geq 0.
\end{eqnarray*}
One deduces that, up to numerical constants, this result is similar to that emphasized in theorem~\ref{theo:result} since for small deviation level, $\alpha (r) = O(\rho_0 r^2 /b)$ whereas $\alpha (r) = O(r/\sqrt{a})$  for large $r$. \vspace{0.1cm}

As we have seen above, the entropic inequality (\ref{eq:mlsi}) entails on the one hand a concentration property for the invariant measure $\mu$ through observables $f\in \L _V(a,b)$. On the other hand and as announced in remark~\ref{rem:poincare}, the Poincar\'e inequality can be seen as a minimal assumption in our study. Hence one may wonder if the Beckner-type inequality (\ref{eq:tetali}), which interpolates between both, provides qualitative concentration estimates through observables in $\L _V(a,b)$. Our second main result, theorem~\ref{theo:poincare}, goes in this way. However, although the estimates we provide below are somewhat similar to that of theorem~\ref{theo:result} in regard of the mixed Gaussian-exponential behaviour, the price to
pay is to lose the good order of magnitude for the Gaussian window since we obtain $r_{\max} = O (\sqrt{b})$ instead of $O(b)$.
\begin{theo}
\label{theo:poincare}
Assume that there exists $p\in (1,2]$ such that the pair $(\mu,\Gamma)$ satisfies a Beckner-type inequality (\ref{eq:tetali}) and denote $\alpha_p$
the optimal constant. Moreover, assume that the observable $f\in \L _V(a,b)$ with the restriction $\alpha _p \leq  2b(p-1)/(3pa)$ and let
$$
r_{\max} := \sqrt{\frac{32bp}{27(p-1)\alpha_p}}
$$
be the size of the Gaussian window. Then the following tail estimates hold: for any $0\leq r \leq r_{\max}$,
\begin{equation}
\label{eq:conc3}
\mu \left( \left\{ x\in \X: f(x) - \mu (f) >r \right\} \right) \leq e^{- \frac{9\alpha_p \, r^2}{32 b}},
\end{equation}
whereas for any $r\geq r_{\max}$,
\begin{equation}
\label{eq:conc4}
\mu \left( \left\{ x\in \X: f(x) - \mu (f) >r \right\} \right) \leq e^{- r \, \sqrt{\frac{3p \alpha_p }{32b(p-1)}}} .
\end{equation}
\end{theo}
\begin{Proof}
The proof is adapted from the method of Aida and Stroock introduced in \cite{aida_stroock}. Assume without loss of generality that $f$ is centered and bounded and for any
$\lambda \in (0,\lambda_0)$, where
$$
\lambda_0 \, := \, \sqrt{\frac{3p \alpha_p}{2b(p-1)}} \, \leq \, \frac{1}{\sqrt{a}},$$
denote once again $Z_\lambda := \int_\X e^{\lambda f} \, d\mu$ and $\mu_\lambda$ the probability measure with density $f_\lambda := e^{\lambda f}/Z_\lambda $ with respect to $\mu$. We have by the Beckner-type inequality (\ref{eq:tetali}) applied to the function $e^{\lambda f/p}$,
\begin{eqnarray*}
Z_\lambda - Z_{\lambda/p} ^p & \leq & \frac{p}{2\alpha_p} \, \cE_\mu (e^{\lambda f/p},e^{\lambda (1-1/p) f}) \\
& = & \frac{p}{2\alpha_p} \, \underset{f(x) >f(y) }{\int \int} e^{\lambda f(x)} \left( 1 - e^{- \lambda (f(x)-f(y))/p}\right) \left( 1- e^{- \lambda (1-1/p) (f(x)-f(y))}\right) \, Q_x (dy) \mu (dx) \\
& \leq & \frac{\lambda^2 (p-1) Z_\lambda}{2p \alpha_p} \, \int_\X  \Gamma (f,f) \, d\mu_\lambda \\
& \leq & \frac{\lambda^2 (p-1) Z_\lambda}{2p \alpha_p} \, \int_\X \left( - a \, \frac{\L V}{V} + b \right) \, d\mu_\lambda \\
& \leq & \frac{\lambda^2 (p-1) Z_\lambda}{2p \alpha_p} \, \left( a I(\mu_\lambda |\mu) +b \right) \\
& \leq & \left( \frac{\lambda }{\lambda_0} \right)^2\, Z_\lambda ,
\end{eqnarray*}
where we used that $f\in \L _V(a,b)$ and lemmas~\ref{lemme:fisher}-\ref{lemme:cle} in the three last lines. Hence rearranging the terms above and iterating the procedure yields for every $n\geq 1$,
\begin{eqnarray*}
Z_\lambda & \leq & \prod_{k=0}^{n-1} \left( \frac{\lambda_0 ^2}{\lambda_0 ^2- \lambda ^2 /p^{2k}} \right) ^{p^k} \, (Z_{\lambda/p^n}) ^{p^n} .
\end{eqnarray*}
Since $f$ is centered, the quantity $Z_{\lambda/p^n} ^{p^n}$ goes to 1 as $n\to \infty$ and from the latter inequality we obtain after taking
logarithm,
\begin{eqnarray*}
\log Z_\lambda & \leq & - \sum_{k=0}^{\infty} p^k \, \log \left( 1 - \frac{(\lambda /\lambda_0) ^2}{p^{2k}} \right) \\
& = & \sum_{k=0}^{\infty} \frac{p^{2k+1}}{p^{2k+1}-1} \times \frac{(\lambda/\lambda_0)^{2(k+1)}}{k+1} \\
& \leq & - \frac{p}{p-1} \, \log \left( 1- \left( \frac{\lambda}{\lambda_0} \right)^2\right).
\end{eqnarray*}
In the last inequality we used the trivial inequality $p^{2k+1} \leq (\frac{p}{p-1}) \, ( p^{2k+1}-1)$ available for any integer $k$
because $p\in (1,2]$. We thus obtain for any $0<\lambda \leq \lambda_0 /2$,
\begin{eqnarray*}
Z_\lambda & \leq & \left( 1+ \frac{\lambda^2}{\lambda_0 ^2 - \lambda^2 }\right) ^{p/(p-1)} \\
& \leq & \exp \left( \frac{p\lambda^2}{(p-1) (\lambda_0 ^2 - \lambda^2)}\right) \\
& \leq & \exp \left( \frac{4p \lambda^2}{3(p-1)\lambda_0 ^2 }\right).
\end{eqnarray*}
Finally using the exponential Chebyshev inequality entails the desired result.
\end{Proof}
\begin{remark}
\rm{The assumption relying $\alpha_p$ to the parameters $a$ and $b$ is a technical detail but cannot be avoided. However it will be always satisfied
as soon as $b$ is taken sufficiently large (or $a$ small enough).}
\end{remark}
\begin{remark}
\rm{As already mentioned, the Beckner-type inequality is stronger than the Poincar\'e inequality, i.e.
$\alpha_p \leq \lambda_1$. However, theorem~\ref{theo:poincare} does not entail a better concentration estimate than that obtained under the
Poincar\'e inequality, except maybe when focusing on the constants depending on $p$ (this is clearly not our interest here). The reason is due to the approach emphasized above which is exactly the same for any $p\in (1,2]$, in contrast to theorem~\ref{theo:result} where the Herbst method is
used.}
\end{remark}

\section{Examples}
\label{sect:ex}

\subsection{Diffusion processes}
Let us apply now theorems~\ref{theo:result} and \ref{theo:poincare} to diffusion processes. Below, the function $U$ is a smooth potential such that $e^{-U}$ is Lebesgue integrable, and denote $\mu$ the Boltzmann probability measure with density $e^{-U}/Z$ with respect to the Lebesgue measure, where $Z$ is the normalization factor. \vspace{0.1cm}

The first example of interest is the so-called Kolmogorov process with generator given for any $C^2$ function $f:\R^d \to \R$ with bounded derivatives by
$$
\L f = \Delta f - <\nabla U , \nabla f> .
$$
One easily checks that $\mu$ is reversible for this process and the carr\'e du champ is $\Gamma (f,f) = \Vert \nabla f\Vert ^2$ where $\Vert \cdot \Vert $ stands for the Euclidean norm in $\R^d$. Hence by Rademacher's theorem, the spaces $\Lip(\X)$ and $\Lip _\Gamma (\X)$ coincide. Moreover the domain $\D (\cE _\mu)$ of the Dirichlet form is $H^1 (\mu)$.

\subsubsection{Ornstein-Uhlenbeck process and the standard Gaussian distribution}
Let us consider first the case of the Ornstein-Uhlenbeck process which has the standard Gaussian distribution as invariant measure. Here the
potential is given by $U(x) = \Vert x\Vert ^2 /2$. By the famous Gross theorem \cite{gross}, the pair $(\mu,\Gamma)$ satisfies the log-Sobolev inequality, i.e. the entropic inequality with (optimal) constant $\rho_0=2$. Hence theorem~\ref{theo:result} will apply for observables in $\L _V(a,b)$ for some good test function $V$. For instance if $f(x) = \Vert x\Vert ^2$, then choose the positive test function $V= e^{cU}$ with $c\in (0,1)$,
i.e. $V$ is at the boundary of non-integrability. Then we have
\begin{eqnarray*}
- \frac{\L V (x)}{V(x)} & = & -c d + c(1-c) \, \Vert x\Vert ^2.
\end{eqnarray*}
Thus with the choice $c=1/2$ we get $f\in \L _V(a,b)$ with $a=16$ and $b=8d$ and by theorem~\ref{theo:result}, for any $0\leq r \leq 8d/3$,
\begin{eqnarray*}
\mu \left( \{ x\in \R^d : \Vert x\Vert ^2 > d+r \} \right) & \leq & e^{- \frac{3r^2}{64d}},
\end{eqnarray*}
which is sharp up to a numerical constant since in this case $\Var_\mu (f) = 2d$. In the exponential regime we get for any $r\geq 8d/3$,
\begin{eqnarray*}
\mu \left( \{ x\in \R^d : \Vert x\Vert ^2 > d+r \} \right) & \leq & e^{- \frac{r}{8}}.
\end{eqnarray*}
Actually, such a behaviour is expected since under $\mu$, the variable $f(x)$ is $\chi ^2 (d)$-distributed, i.e. distributed as a chi-squared random variable with $d$ degrees of freedom. \vspace{0.1cm}

Let us come back to the example announced in remark \ref{rem:regimes}. As observed above, our concentration result is convenient as soon as $f$ is close to realize the equality in the Lyapunov condition (\ref{eq:lyap}). However what happens for a function $g$ such that $\Vert \nabla g \Vert \ll \Vert \nabla f \Vert $ at
infinity ?
For instance if $f(x)= \Vert x\Vert ^2$ as above then how concentrates the invariant measure through the observable $g(x) = \Vert x\Vert ^{3/2}$ ? Let us investigate this point in detail now. Assume that there exists an observable $g$ satisfying for any $\varepsilon >0$,
\begin{eqnarray*}
\Vert \nabla g \Vert ^2 & \leq & \varepsilon \Vert \nabla f\Vert ^2 + \frac{1}{\varepsilon } \\
& \leq & - a_{\varepsilon} \, \frac{\L V}{V} + b_\varepsilon ,
\end{eqnarray*}
where $a_\varepsilon := a \varepsilon$ and $b_\varepsilon := b\varepsilon + 1/\varepsilon$. Using the argument given in the proof of lemma~\ref{lemme:cle}, we have for any $\lambda>0$,
\begin{eqnarray*}
I(\mu _\lambda |\mu) & \leq & \inf _{\varepsilon >0} \, \frac{\lambda ^2}{4} \, \left( a_\varepsilon  I(\mu _\lambda |\mu) + b_\varepsilon  \right) \\
& = & \frac{\lambda ^2}{2} \, \sqrt{a I(\mu _\lambda |\mu) + b},
\end{eqnarray*}
where in the definition of $I(\mu _\lambda |\mu)$ we replaced $f$ by $g$. Hence we obtain
\begin{eqnarray*}
I(\mu _\lambda |\mu) & \leq & \frac{\lambda ^2}{2} \, \left( \frac{a\lambda^2}{2} + \sqrt{b} \right) .
\end{eqnarray*}
Now the same argument as in the proof of theorem~\ref{theo:result} together with the latter inequality entail
\begin{eqnarray*}
\frac{d}{d\lambda } L_\lambda & \leq & \inf _{\varepsilon >0} \, \frac{1}{\rho_0} \, \left( a_\varepsilon  I(\mu _\lambda |\mu) + b_\varepsilon
\right) \\
& = & \frac{2}{\rho_0} \, \sqrt{a I(\mu _\lambda |\mu) + b} \\
& \leq & \frac{2}{\rho_0} \, \left( \frac{a\lambda^2}{2} + \sqrt{b} \right) ,
\end{eqnarray*}
since $\sqrt{x(x+y) + y^2} \leq x+y$ for any $x,y\geq 0$. Hence we get for any $\lambda>0$,
\begin{eqnarray*}
\log \int_\X e^{\lambda g } \, d\mu & \leq & \lambda \mu (g) + \frac{2}{\rho_0} \, \left( \frac{a\lambda^4}{6} + \lambda ^2 \sqrt{b} \right).
\end{eqnarray*}
Finally, a bit more of analysis shows that the decay in the concentration estimate is of order $e^{-c r^{4/3}}$ for large deviation level $r$, which is of the good order of magnitude when choosing for instance $g(x)$ to be proportional to $\Vert x\Vert ^{3/2}$. \vspace{0.1cm}

Now take the observable $f$ as a quadratic form on $\R^d$, i.e. there exists a positive definite symmetric matrix $A=(a_{i,j})_{i,j=1,\ldots,d}$ of size $d$ such $f(x) = <A x,x>$, $x\in \R^d$. Then in the Gaussian regime we should obtain the variance of $f$ in the denominator,
$$
\Var_\mu (f)  \,  =\, 2 \, \sum_{i,j=1} ^d a_{i,j} ^2 .
$$
However theorem~\ref{theo:result} does not yield such a result, unfortunately. For instance using the same test function $V$ as before would entail that $f\in \L _V(a,b)$ with $a=16\, \Vert A\Vert _{\mathrm{op}} ^2$ and $b= 8d\, \Vert A\Vert _{\mathrm{op}} ^2$, where $\Vert A\Vert _{\mathrm{op}}$ is the (Euclidean) operator norm of the matrix $A$, i.e. its spectral radius. With this choice of parameters the inequality $\Var_\mu (f) \leq b/\rho_0$ is too weak to provide a reasonable variance estimate since $b$ behaves badly in terms of dimension. To circumvent this difficulty, we choose $V= e^{cf}$ with $c=1/(4 \Vert A\Vert _{\mathrm{op}})$, so that $V$ is positive and still integrable with respect to the Gaussian measure $\mu$. We thus obtain that $f\in \L _V(a,b)$ with the same $a$ as before, but now with the improved variance estimate
$$
b = 2ac \,  \rm{trace}(A) = 8 \, \rm{trace}(A) \, \Vert A\Vert _{\mathrm{op}},
$$
where $\rm{trace}(A)$ is the trace of the matrix $A$. Therefore applying theorem~\ref{theo:result} entails a tail estimate whose behaviour has been improved with respect to the dimension. This example emphasizes the inherent importance of the choice of the function $V$ in the condition $\L_V(a,b)$. See also for instance \cite{hanson, latala_chaos} for some nice studies on the concentration properties of Gaussian-like quadratic forms and Gaussian chaoses.

\subsubsection{Kolmogorov process and the Boltzmann invariant measure} This class is a natural generalization of the Ornstein-Uhlenbeck process. To begin, assume that the measure $\mu$ is spherically log-concave, i.e. there exists a $C^2$ function $\phi : \R \to \R$ convex and non-decreasing such that $U(x) = \phi(\Vert x\Vert )$ for any $x\in \R^d$. By a famous result of Bobkov \cite{bobkov}, the dynamics $(\mu, \Gamma)$ satisfy (at least) a Poincar\'e inequality. Let us consider the potential $U$ as an observable and also the positive test function $V = e^{cU}$, which belongs to $L^1(\mu)$ for any $c\in (0,1)$ since $\phi$ is convex. Assume that there exists $r>0$ and $M= M(r) \in (0,1-c)$ such that
$$
\frac{\Delta U (x) }{\Vert \nabla U (x) \Vert^2 } \, = \, \frac{(d-1) \, \phi '(\Vert x\Vert) + \Vert x\Vert \, \phi ''(\Vert x\Vert)}{\Vert x\Vert \, \phi ' (\Vert x\Vert) ^2} \, \leq \, M , \quad \Vert x\Vert \geq r .
$$
Since we have
\begin{eqnarray*}
- \frac{\L V (x)}{V(x)} & = & - c \, \Delta U (x) + c(1-c) \, \Vert \nabla U(x) \Vert ^2 \\
& = & -c (d-1) \, \frac{\phi '(\Vert x\Vert)}{\Vert x\Vert} - c \, \phi ''( \Vert x\Vert) +c(1-c) \, \phi ' (\Vert x\Vert) ^2 ,
\end{eqnarray*}
one deduces that $U$ belongs to the class $\L _V (a,b)$ with $a= 1/c(1-c-M)$, the parameter $b=b(r)$ being chosen conveniently on $B_r$, the centered ball of radius $r$ in $\R^d$, i.e. $b$ is the maximum between $3a\lambda_1$ and   
$$
\left \Vert \left( a \, \frac{\L V}{V} + \Vert \nabla U\Vert ^2 \right) \, 1_{ B_r } \right \Vert _{L^\infty (\mu)} .
$$
As a result, one can apply theorem~\ref{theo:poincare} to obtain a Gaussian-exponential concentration estimate through the observable $U$. See also the recent work of Bobkov and Madiman \cite{bob_madiman} where a somewhat similar tail estimate is established via a different approach. \vspace{0.1cm}

Actually, spherically log-concave probability measures include the case of a potential $U$ such that $U(x) = \Vert x\Vert ^\beta$ with $\beta \geq 1$. Since the case $\beta = 2$ has already been considered, three different situations arise:

$(i)$ the case $\beta =1$, for which only the Poincar\'e inequality is satisfied.

$(ii)$ the case $\beta \in (1,2)$: the standard Beckner inequality holds, cf. \cite{latala}.

$(iii)$ the case $\beta > 2$: using Wang's criterion \cite{wang}, the log-Sobolev inequality is then verified. \vspace{0.1cm}

In these three cases, one may choose the following parameters:
$$
c := \frac{1}{2}, \quad M := \frac{1}{4} , \quad a := 8 \quad \textrm{ and } \quad r:= 4 ^{1/\beta} \, \left( \frac{d+\beta -2}{\beta} \right) ^{1/\beta} ,
$$
provided the restriction $d+\beta -2 >0$ holds. Finally, if $\beta >2$ then the parameter $b$ can be easily chosen, in contrast to the case $\beta \in [1,2)$ for which $U$ is not $C^2$ at 0. Therefore, to obtain a convenient constant $b$ in this non-smooth situation, one can choose a test function $V=e^{c\tilde{U}}$ where $\tilde{U}$ is $C^2$ on $\R^d$ and $U=\tilde{U}$ outside the ball $B_r$. Then the proof above remains valid with $\tilde{U}$ instead of $U$ and an easy perturbation argument entails the standard Beckner inequality (or the Poincar\'e inequality in the case $\beta =1$) for the Boltzmann probability measure defined with respect to the potential $\tilde{U}$. \vspace{0.1cm}

One may also extend the result to the non symmetrically invariant case, when for example a logarithmic Sobolev inequality holds. Let us assume for example that the potential $U$ is such that its Hessian matrix, denoted $\mbox{Hess}\, U$, is lower bounded and that the following Lyapunov condition holds:
$$
\L V (x) \leq \left( - c_1 \, \Vert x\Vert ^2 + c_2 \right) \, V (x), \quad x\in \R^d,
$$
where $c_1,c_2 >0$ and $V$ is a $C^2$ positive test function. Then a logarithmic Sobolev inequality holds, cf. \cite{cgw}, and one can apply theorem \ref{theo:result} for observables $f$ such that the norm of $\nabla f (x) $ is at most $\Vert x \Vert $ at infinity, since in this case $f\in \L _V (a,b)$ with $a = 1/c_1$ and $b=c_2 /c_1$. For instance, the Lyapunov condition above will be verified if at least one of the two conditions below is satisfied: there exist $\alpha \in (0,1)$ and $\beta >0$ such that for sufficiently large $x$,
$$
(1-\alpha) \, \Vert \nabla U (x) \Vert ^2 - \Delta U(x) \, \geq \,  \beta \Vert x\Vert ^2 \quad \mbox{  or  } \quad < x, \nabla U (x) > \, \geq \, \beta \Vert x\Vert ^2.
$$

\subsubsection{Log-Sobolev inequality for modified dynamics}
Our last example concerns the case where a log-Sobolev inequality holds for a slightly modified dynamics, but with the same Boltzmann invariant measure. In a sense, it corresponds to a modified (or weighted) log-Sobolev inequality. Consider the process with the following generator:
$$
\L^{\sigma^2} f := \sigma ^2 \, \Delta f + <\nabla (\sigma^2) - \sigma^2 \nabla U , \nabla f> ,
$$
where $\sigma$ is some measurable and locally bounded function from $\R^d$ to $\R$. Once again the measure $\mu$ is reversible for this process,
but the notable difference relies on the weight $\sigma^2$ in the carr\'e du champ, i.e. $\Gamma^{\sigma^2} (f,f) := \sigma ^2 \, \Vert \nabla f\Vert ^2$, so that a Lipschitz function $f$ may have an unbounded carr\'e du champ $\Gamma^{\sigma^2}(f,f)$, in contrast to the Kolmogorov example studied above. In particular, the domain of the Dirichlet form is a weighted $H^1$ space, i.e.
$$
\D (\cE _\mu) := \left \{ f\in L^2(\mu): \int_{\R^d} \sigma ^2 \, \Vert \nabla f \Vert ^2 \, d\mu <\infty \right \}.
$$
We will focus mainly here on the simple case $U(x) := \Vert x\Vert ^\alpha$ for $1<\alpha<2$, so that the standard Beckner inequality (thus the Poincar\'e inequality) holds for the classical dynamics $(\mu, \Gamma)$, but not a log-Sobolev inequality. On the one hand, according to a result of Latala and Oleszkiewicz \cite{latala}, the measure $\mu$ concentrate like $e^{-r ^\alpha} $ for large deviation level $r$ through Lipschitz observables. On the other hand, the following weighted log-Sobolev inequality holds:
$$
\mbox{Ent}_\mu (f^2) \leq C \, \int_{\R^d} \left( 1+\Vert x\Vert ^{2-\alpha}\right) \, \Vert \nabla f(x) \Vert ^2 \, \mu (dx),
$$
where $C >0$ is some constant depending on dimension $d$, cf. \cite{cgw2}. Letting the weight function $\sigma(x) ^2:= 1+\Vert x\Vert ^{2-\alpha}$, one observes that this weighted inequality rewrites as the log-Sobolev inequality for the new dynamics $(\mu,\Gamma ^{\sigma^2})$. Choosing the positive test function $V(x) = e^{c\Vert x\Vert ^\alpha}$, which belongs to $L^1(\mu)$ for any $c\in (0,1)$, we have for any $x$ outside a neighborhood of 0,
\begin{eqnarray*}
- \, \frac{\L^{\sigma^2} V(x)}{V(x)} & = & -c\alpha (d+\alpha-2) \, (1+\Vert x\Vert ^{2-\alpha}) \Vert x\Vert ^{\alpha-2} + \alpha^2 c(1-c) (1+\Vert x\Vert ^{2-\alpha}) \Vert x\Vert ^{2(\alpha-1)} \\
& &  - c \alpha (2-\alpha),
\end{eqnarray*}
which behaves like $\alpha^2 c(1-c) \Vert x\Vert ^\alpha$ at infinity. Hence using the same reasoning as in the case of the Kolmogorov process above, one deduces that observables $f$ having a gradient $\Vert \nabla f(x) \Vert$ controlled by $\Vert x\Vert ^{\alpha -1}$ for large $x$ satisfy theorem~\ref{theo:result} (the observable $f(x)=\Vert x\Vert ^\alpha$ belongs to this class, as expected according to \cite{latala}). We point out that our results might be made more precise by following the approach provided in remark~\ref{rem:regimes}. To that aim, one has to consider the functional inequality $I_\mu (a)$ involved in Latala and Oleszkiewicz's work \cite{latala}, which is more general than the standard Beckner inequality emphasized above. \vspace{0.1cm}

In fact using the modified dynamics $(\mu, \Gamma^{\sigma^2})$, one can even consider interesting cases for which even the Poincar\'e inequality does not hold for the original dynamics $(\mu, \Gamma)$. For instance consider the generalized Cauchy measure $\mu$ with density proportional to $(1+\Vert x\Vert ^2)^{-\beta}$, where the condition $\beta >d/2$ holds to ensure integrability at infinity. Such a measure satisfies both a weighted Poincar\'e inequality:
\begin{equation}
\label{eq:poincare_cauchy}
\Var _\mu (f) \leq C \, \int _{\R ^d} \left( 1+\Vert x\Vert ^2 \right) \, \Vert \nabla f(x) \Vert ^2 \, \mu (dx),
\end{equation}
with the slight restriction $\beta \geq d$, cf. \cite{bob_ledoux2}, and also a weighted log-Sobolev inequality according to \cite{cgw2}:
\begin{equation}
\label{eq:lsi_cauchy}
\Ent _\mu (f^2) \leq  \tilde{C} \, \int _{\R ^d} \left( 1+\Vert x\Vert ^2 \right) \, \log \left( e+\Vert x\Vert ^2 \right) \, \Vert \nabla f(x) \Vert ^2 \, \mu (dx) ,
\end{equation}
without further restriction on $\beta$. Here $C$ and $\tilde{C}$ are some positive constants depending on $\beta$ and $d$. Letting the weights
$$
\sigma_1 (x) ^2 := 1+\Vert x\Vert ^2 \quad \mbox{  and  } \quad \sigma _2(x) ^2 := \left( 1+\Vert x\Vert ^2 \right) \, \log \left( e+\Vert x\Vert ^2 \right) ,
$$
then the weighted inequalities (\ref{eq:poincare_cauchy}) and (\ref{eq:lsi_cauchy}) rewrite as the Poincar\'e and the log-Sobolev inequality for the dynamics $(\mu, \Gamma ^{\sigma_1 ^2})$ and $(\mu, \Gamma ^{\sigma_2 ^2})$, respectively. Now let $V(x) = \Vert x\Vert ^k$ for some $0 < k < 2\beta -d $, so that the positive test function $V$ lies in $L^1(\mu)$. Then we have for any $x$ outside a neighborhood of 0,
\begin{eqnarray*}
- \, \frac{\L^{\sigma_1 ^2} V(x)}{V(x)} & = & - k (d+k-2) \, \frac{1+ \Vert x\Vert ^2}{\Vert x\Vert ^2} + 2k(\beta -1) ,
\end{eqnarray*}
and also
\begin{eqnarray*}
- \, \frac{\L^{\sigma_2 ^2} V(x)}{V(x)} & = & - k (d+k-2) \, \frac{1+ \Vert x\Vert ^2}{\Vert x\Vert ^2} \, \log \left( e+\Vert x\Vert ^2 \right) - \frac{2k\left( 1+ \Vert x\Vert ^2 \right)}{e+ \Vert x\Vert ^2} \\
& & + 2k(\beta -1) \, \log \left( e+\Vert x\Vert ^2 \right) .
\end{eqnarray*}
The first quantity is constant at infinity whereas the second one is of order $k(2\beta -d - k )  \, \log \left( e+\Vert x\Vert ^2 \right)$ for large $\Vert x\Vert$. Then we obtain by theorems~\ref{theo:result} and \ref{theo:poincare} an exponential concentration estimate for large deviation level $r$ through observables $f$ having their gradient $\Vert \nabla f(x) \Vert$ dominated in both cases by $1/ \Vert x\Vert $ for large $x$. Hence the previous example shows that weighted Poincar\'e and log-Sobolev inequalities (of course with a different weight) can lead to somewhat similar concentration estimates. Note that the function $f(x)= \log \left( \Vert x\Vert \right)$ belongs to this class of observables, leading to the well-known heavy tail phenomenon satisfied by Cauchy-type measures, cf. \cite{bob_ledoux2}. Finally, we mention that one can take profit of remark~\ref{rem:regimes} to get intermediate concentration regime for observables not saturating the Lyapunov condition.

\subsection{Birth-death processes}
Let us begin the study of jump processes by considering a simple but however non trivial example, namely birth-death processes. Here $(X_t)_{t\geq 0}$ is a Markov process on the state space $\N := \{ 0,1,2, \ldots \}$ endowed with the classical metric $d(x,y) = \vert x-y\vert$, $x,y\in \N$. The transition probabilities are given by
\[
\P_x (X_t =y) = %
\begin{cases}
  \lambda_x t + o(t) & \text{if $y=x+1$}, \\
  \nu_x t + o(t) & \text{if $y=x-1$}, \\
  1- (\lambda_x + \nu_x) t + o(t) & \text{if $y=x$},
\end{cases}
\]
where $\lim_{t\to0}t^{-1}o(t)=0$. The transition rates $\lambda$ and $\nu$ are respectively called the birth and death rates and satisfy to $\lambda >0$ on $\N$ and $\nu >0$ on $\N ^* := \{ 1,2,\ldots \}$ and $\nu_0 = 0$, so that the process is irreducible. Although we assume that the stability condition (\ref{eq:stab}), which rewrites as
$$
\lambda _x + \nu _x < \infty, \quad x\in \N ,
$$
is satisfied, the generator might be unbounded in the sense of (\ref{eq:unbounded}), i.e.
\begin{eqnarray*}
\sup _{x\in \N} \, \lambda_x + \nu_x = \infty .
\end{eqnarray*}
The process is positive recurrent and non-explosive when the rates satisfy to
\[
\sum_{x=1}^\infty \frac{\lambda_0 \lambda_1 \cdots \lambda_{x-1}}{\nu_1 \nu_2 \cdots \nu_x} <\infty
\quad\text{and}\quad
\sum_{x=1}^\infty %
\left(\frac{1}{\lambda_x}+\frac{\nu_x}{\lambda_x \lambda_{x-1}}%
  +\cdots+\frac{\nu_x\cdots\nu_1}{\lambda_x\cdots \lambda_1 \lambda_0}\right) = \infty,
\]
respectively. In this case the detailed balance condition (\ref{eq:revers}) rewrites as
$$
\lambda _x \, \mu(\{x\}) = \nu_{x+1} \, \mu(\{ x+1\}) , \quad x\in \N ,
$$
where $\mu$ is the unique stationary distribution of the process given by
\begin{equation}\label{eq:invariant}
\mu (\{x\}) = \mu (\{0\}) \prod_{y=1}^x \frac{\lambda _{y-1}}{\nu_y} ,\ x\in\N ,
\end{equation}
$\mu(\{ 0\})$ being the normalization constant. In the situations of interest, the death rate $\nu$ has to be bigger than $\lambda$ to ensure such criteria. \vspace{0.1cm}

For any function $f:\N \to \R$, the generator $\L$ of the process is given by
\begin{eqnarray*}
\L f(x) & = & \lambda_x \, \left( f(x+1) -f(x)\right) + \nu_x \, \left( f(x-1) -f(x)\right), \quad x\in \N,
\end{eqnarray*}
and the carr\'e du champ is
\begin{eqnarray*}
\Gamma (f,f)(x) & = & \frac{1}{2}\,\left\{ \lambda_x \, \left( f(x+1)-f(x) \right) ^2 + \nu_x \left( f(x-1)-f(x) \right) ^2 \right\} , \quad x\in \N.
\end{eqnarray*}
In particular, the Dirichlet form is given by
$$
\cE _\mu (f,g) := \sum_{x\in \N} \lambda _x \, \left( f(x+1)-f(x) \right) \, \left( g(x+1)-g(x) \right) \, \mu (\{x\}) ,
$$
where $f,g$ belong to the space $\D (\cE_ \mu)$ of functions $u: \N \to \R$ such that $\cE _\mu (u,u)$ is finite. \vspace{0.1cm}

On the one hand Joulin \cite{joulin} gives, under some convenient ergodic assumptions, concentration estimates of Poisson-type through observables belonging to the space $\Lip(\N)$. In particular, his proof requires the inclusion $\Lip(\N) \subset \Lip_\Gamma (\N)$, enforcing the rates $\lambda$ and $\nu$ to be bounded. On the other hand, when we apply Ollivier's result \cite{ollivier} to birth-death processes, his Gaussian-exponential concentration property is available for observables $f\in \Lip(\N)$ such that $\Gamma(f,f) \in \Lip(\N)$. It induces that $\lambda, \nu\in \Lip(\N)$, extending Joulin's result from bounded to (at most) linear rates. As announced in the introduction through the study of a specific example, theorems~\ref{theo:result} and \ref{theo:poincare} entail Gaussian-exponential concentration estimates beyond these cases since the carr\'e du champ $\Gamma(f,f)$ is allowed to have a growth comparable to that of $\nu$. \vspace{0.1cm}

First let us provide some basic conditions which ensure an entropic or Poincar\'e inequality. The following necessary (but not sufficient) condition is due to Caputo, DaiPra and Posta \cite{caputo} and has been recently recovered by Chafa\"i and Joulin \cite{cha_jou} by using a semigroup approach: if $\lambda$ is non-increasing and $\nu$ is non-decreasing and there exists $\alpha >0$ such that
\begin{eqnarray}
\label{eq:daipra}
\inf_{x\in \N} \, \lambda _x - \lambda_{x+1} + \nu_{x+1}-\nu_x \geq \alpha,
\end{eqnarray}
then the entropic inequality (\ref{eq:mlsi}) is satisfied with constant $\alpha$, or equivalently $\rho_0 \geq \alpha$. Such an assumption exhibits very asymmetric rates. More precisely, it enforces the rates $\lambda$ and $\nu$ to be bounded and super-linear, respectively, excluding some interesting cases which can be however considered for the Poincar\'e inequality. Indeed, Miclo \cite{miclo} states that the spectral gap $\lambda_1$ is positive if and only if
\begin{eqnarray}
\label{eq:miclo}
\delta :=  \sup _{x\geq 1} \, \sum_{k=0} ^{x-1} \frac{1}{\lambda _k \mu (\{k\})} \, \sum_{l\geq x} \mu (\{ l \}) & <& \infty,
\end{eqnarray}
and in this case we have $1/\delta \geq \lambda_1 \geq 1/4\delta$, i.e. $\lambda_1$ is of order $1/\delta$. In particular, in contrast to the entropic inequality, one may find examples satisfying the Poincar\'e inequality with an unbounded birth rate $\lambda$. Now assume that the positive test function $V(x):= \kappa ^{x}$ is in $L^1(\mu)$ for some constant $\kappa >1$ depending on $\lambda,\nu$. Then an observable $f$ belongs to $\L_V (a,b)$ if and only if
\begin{eqnarray*}
\Gamma (f,f) & \leq & \frac{a (\kappa -1)}{\kappa} \, \left( \nu - \kappa \, \lambda \right) + b ,
\end{eqnarray*}
showing that on a large scale the behaviour of $\Gamma (f,f)$ is controlled by the growth of the death rate $\nu$. To compare with the results of Joulin and Ollivier mentioned previously, assume that the observable $f\in \Lip (\N)$. Then two extreme situations may appear when the death rate $\nu$ is unbounded:

$(i)$ a small birth rate $\lambda$, i.e. $\lambda$ is bounded. In this case one may choose the following parameters to ensure that $f\in \L_V (a,b)$:
$$
a := \frac{\kappa}{2(\kappa-1)} \quad \mbox{and} \quad b := \frac{(1+\kappa) \Vert \lambda \Vert _{L^\infty (\mu)}}{2} .
$$

$(ii)$ a birth rate $\lambda$ of the order of $\nu$. Let $x_0 \in \N ^*$ and assume that $\lambda _x \leq c \nu_x$ for all $x\geq x_0$, where $c\in (0,1)$ is some parameter. Then in order to get $f\in \L_V (a,b)$, one can choose for $\kappa \in (1,1/c)$,
$$
a:= \frac{\kappa (1+c)}{2(1-c\kappa)(\kappa -1)},
$$
$$
\mbox{and} \quad \quad b := \left \Vert \left( \frac{\lambda +\nu}{2} +a \, \frac{\L V}{V } \right) \, 1 _{[0,x_0 ]} \right \Vert  _{L^\infty (\mu)} = \, \frac{1+\kappa}{2(1-c\kappa )} \, \left \Vert (\lambda - c \nu ) \, 1 _{[0,x_0 ]}\right \Vert _{L^\infty (\mu)} .
$$
For instance the choice of $\kappa := 1/\sqrt{c}$ entails the integrability of $V$ and then we obtain
$$
a = \frac{1+\sqrt{c}}{2(1-\sqrt{c})} \quad  \quad \mbox{and} \quad \quad b = \frac{1+\sqrt{c}}{2\sqrt{c} (1-\sqrt{c})}  \, \left \Vert (\lambda - c \nu ) \, 1 _{[0,x_0 ]}\right \Vert _{L^\infty (\mu)} .
$$ 
Of course $b$ has also to be at least $3a\lambda _1$ if only the Poincar\'e inequality is satisfied. In both cases $(i)$ and $(ii)$ there exist plenty of examples satisfying Poincar\'e inequality and thus theorem~\ref{theo:poincare}, whereas only the case $(i)$ may satisfy the entropic inequality and so theorem~\ref{theo:result}. In particular, the example emphasized in the introduction is a prototype of such a situation. Indeed denote $\mu_p$ the geometric distribution on $\N$ of parameter $p \in (0,1)$, i.e. $\mu_p (\{ x \} ) := (1-p) p^x$, $x\in \N$, and let us consider the carr\'e du champ
$$
\Gamma ^{(n)} (f,f) (x) := \frac{1}{2} \, \left \{ p(x+1)^n \, (f(x+1) - f(x)) ^2 + x^n 1_{\{x \neq 0\}} \, (f(x-1) - f(x)) ^2\right\},
$$
where $x\in \N$ and $n\in \N$ is some fixed parameter. The measure $\mu_p$ is reversible with respect to these dynamics and by comparing the underlying Dirichlet form with that given in the case $n=0$ (for which $\lambda_1 ^{(1)}= (1-\sqrt{p})^2$, cf. \cite{chen1, joulin0}), the dynamics $(\mu_p,\Gamma ^{(n)})$ satisfies a Poincar\'e inequality with optimal constant $\lambda_1 ^{(n)} \geq (1-\sqrt{p})^2$ (Miclo's result \eqref{eq:miclo} would only entail $\lambda_1 ^{(n)} \geq (1-p)^2 /4$). Hence for sufficiently large integer $x_0$ one can pick $c= (1+p)^2 /4 \in (0,1)$, $\kappa = 2/(1+p)$ and for some integer $x_1 \in (0,x_0)$,
$$
a = \frac{3+p}{2(1-p)} \quad \quad \mbox{ and } \quad \quad b = \max \left \{ \frac{3+p}{4(1-p)^2} \, \left \vert 4 p (x_1 +1)^n - (1+p)^2 \, x_1 ^n \right \vert , 3a \lambda_1 ^{(n)} \right \}.
$$ 

To achieve the birth-death example, let us focus our attention on a model which mimics the diffusion case, namely ultra log-concave distributions on $\N$, see for instance \cite{johnson, caputo}. We say that a probability measure $\mu$ on $\N$ is ultra log-concave (resp. log-concave) if it satisfies for any $x\in \N ^*$,
$$
x\, \mu (\{x\}) ^2 \geq (x+1) \, \mu (\{x+1\}) \, \mu (\{x-1\}) \quad \mbox{(resp.} \quad \mu (\{x\}) ^2 \geq \mu (\{x+1\}) \, \mu (\{x-1\}) \mbox{)}.
$$
For instance the Poisson distribution is ultra log concave whereas the geometric measure is only log-concave. Assume that the measure $\mu$ has density $e^{-U}/Z$ with respect to the counting measure on $\N$, where $U$ is some nice function and $Z$ is the normalization constant. Denote $\Delta U$ the discrete Laplacian of the potential $U$, i.e.
$$
\Delta U(x) := U(x+1) - 2 U(x) + U(x-1), \quad x\in \N ^* .
$$
Then $\mu$ is ultra log-concave (resp. log-concave) if and only if $\Delta U(x)\geq \log (1+1/x)$ for any integer $x\in \N ^*$ (resp. $\Delta U$ is non-negative). \vspace{0.1cm} \\
From a dynamical point of view, measure $\mu$ is the stationary distribution of the birth-death process with rates
$$
\lambda _x = 1 \quad \mbox{and} \quad \nu_x = e^{U(x)-U(x-1)} \, 1_{\{ x\neq 0\}} , \quad x\in \N .
$$
Then under the ultra log-concavity assumption, we have for any integer $x\geq 2$,
\begin{eqnarray*}
\lambda _x - \lambda _{x+1} + \nu_{x+1} - \nu_x & = & \left( e^{\Delta U(x)} -1\right) \, e^{\sum_{k=1}^{x-1} \Delta U(k) +U(1)-U(0)} \\
& \geq & \prod _{k=1} ^{x-1} \left( 1+\frac{1}{k}\right) \, \frac{e^{U(1)-U(0)}}{x} \\
& = & e^{U(1)-U(0)} \\
& = & \nu_1,
\end{eqnarray*}
so that (\ref{eq:daipra}) is satisfied with $\alpha = \nu_1$ (the cases $x\in \{ 0,1\}$ being straightforward). Thus the entropic inequality holds with constant $\rho_0 \geq \nu_1$. In particular, the super-linearity of the death rate $\nu$ entails that the potential $U$ has a growth at infinity at least $x\log (x)$, showing that the tail behaviour of $\mu$ can be compared to that of a Poisson distribution. Finally, note that the log-concavity assumption only is not sufficient to ensure an entropic inequality since one obtains in this case
$$
\inf_{x\in \N} \, \lambda _x - \lambda _{x+1} + \nu_{x+1} - \nu_x \geq 0.
$$
However, as in the diffusion case, one may find examples of log-concave distributions on $\N$ satisfying the Poincar\'e inequality by using Miclo's condition (\ref{eq:miclo}), which simply rewrites as
$$
\sup_{x\in \N ^*} \, \underset{0\leq k\leq x-1 < l}{\sum} e^{U(k)-U(l)} \, < \, \infty.
$$

\subsection{Glauber dynamics for unbounded particles}
We consider the situation where $\X$ is the unbounded configuration space $\N ^{\Lambda}$, where $\Lambda$ is a bounded subset of $\Z^d$. For each site $x\in \Lambda$, denote $\eta_x$ the number of particles located at $x$. Given a bounded function $\lambda : \Z ^d \to [0,\infty)$, let $\pi$ be the Poisson measure on $\N ^\Lambda$ with parameter $\lambda$, that is to say
$$
\pi (\{\eta\}) \, = \, \prod _{x\in \Lambda}  e^{-\lambda (x)} \, \frac{\lambda (x) ^{\eta _x}}{\eta _x ! } , \quad \eta \in \N ^\Lambda .
$$
We equip $\N^{\Lambda}$ with the total variation distance which counts the number of different particles. In other words, if $\eta$ and $\bar{\eta}$ are two configurations in $\N ^\Lambda$, then the total variation distance is given by
\begin{eqnarray*}
  d(\eta,\bar{\eta}) & := & \sum_{x\in \Lambda} \vert \eta_x - \bar{\eta} _x\vert .
\end{eqnarray*}
Our definition is a straightforward generalization of the classical notion of total variation distance between probability measures, since it coincides with the usual definition when the configurations are normalized by their total masses. For any $f: \N ^\Lambda \to \R$,
the discrete gradient operators are defined by
\begin{eqnarray*}
D_x ^+ f(\eta) & := & f(\eta +\delta _x)-f(\eta) , \quad D_x ^- f(\eta) \; \; := \; \; f(\eta -\delta _x)-f(\eta), \quad \eta \in \N ^\Lambda ,
\end{eqnarray*}
where $\delta_x$ is the Dirac mass at point $x\in \Lambda$ and by convention $D_x ^- f(\emptyset) := 0$. Note that by \cite{joulin_rub}, a given function $f$ belongs to the space $\Lip (\N^\Lambda)$ if and only if
$$
\sup _{(\eta ,x) \in \N^\Lambda \times \Lambda} \vert D_x ^+ f(\eta ) \vert <\infty.
$$
Now let $\phi : \Z ^d \to [0,\infty)$ be an even function, null at the origin and summable on $\Z^d$, i.e. satisfying $\sum_{x\in \Z ^d} \phi (x) < \infty$. We define the Hamiltonian $H :\N ^\Lambda \to \R_+$ as
\begin{eqnarray*}
H (\eta ) & := & \frac{1}{2} \, \sum_{x,y\in \Lambda} \, \phi (x-y) \, \eta_x \, \eta _y .
\end{eqnarray*}
Then the Gibbs measure $\mu$ at inverse temperature $\beta >0$ is the probability measure on $\N ^\Lambda$ given by
\begin{eqnarray*}
\mu (\{\eta\}) & = & \frac{1}{Z} \, e^{-\beta \, H(\eta)} \, \pi (\{\eta\}),
\end{eqnarray*}
where $Z$ is the normalization constant. As observed below, our study is based on the configuration space $\N ^\Lambda$ since our model exhibits free boundary condition, that is to say $\Lambda$ is, in some sense, disconnected from the lattice $\Z^d$. However the aforementioned model might be extended outside $\Lambda$ by introducing an appropriate boundary condition. \vspace{0.1cm}

Now, let us introduce the Glauber dynamics associated to the Gibbs measure above, which can be seen as a spatial birth-death process, cf. \cite{preston}. If $\eta$ is the configuration of the system at time $t$, then a particle appears or disappears at site $x\in \Lambda$ with rates $\lambda (x) \, e^{-\beta D_x ^+ H(\eta)} dt$ and $dt$, respectively. In particular, the case $H=0$ corresponds to the non-interacting case. The generator $\L$ is thus of birth-death type and defined for any function $f:\N ^\Lambda \to \R $ by
\begin{eqnarray*}
\L (\eta) & := & \sum_{x\in \Lambda} \left( c^- (\eta, x) \, D_x ^- f(\eta ) + c^+ (\eta, x) \, D_x ^+ f(\eta) \right) , \quad \eta \in \N ^\Lambda ,
\end{eqnarray*}
where the rates of the dynamics $c^+$ and $c^-$ are given by
\[
\left \{
\begin{array}{lll}
c^+ (\eta,x) & = & \lambda (x) \, e^{-\beta D_x ^+ H(\eta)} \, = \, \lambda (x) \, e^{-\beta \sum_{y\in \Lambda} \phi(x-y) \, \eta_y }; \\
c^- (\eta,x) & = & \eta_x .
\end{array}
\right.
\]
In particular, the stability condition (\ref{eq:stab}) is clearly satisfied since $\Lambda$ is finite and moreover, according to the detailed balance condition (\ref{eq:revers}) which in our context rewrites as
\begin{eqnarray*}
c^{\pm} (\eta,x) \,\mu (\{\eta\}) & = & c^{\mp} (\eta \pm \delta _x ,x) \, \mu(\{\eta \pm \delta _x \}) , \quad \eta_x >0, \quad (\eta , x) \in \N ^\Lambda \times \Lambda ,
\end{eqnarray*}
the Gibbs measure $\mu$ is reversible for these dynamics. Finally, the carr\'e du champ of an observable $f$ is given by
$$
\Gamma (f,f) (\eta) \, = \, \frac{1}{2} \, \sum_{x\in \Lambda} \left\{ c^- (\eta, x) \, \vert D_x ^- f(\eta )\vert ^2 + c^+ (\eta, x) \, \vert D_x ^+ f(\eta) \vert ^2 \right\} , \quad \eta \in \N ^\Lambda.
$$

Recently, the problem of finding the speed of convergence to equilibrium of this model has been addressed in several articles, cf. for instance \cite{bertini} or \cite{wu2} for a spectral method (i.e. related to Poincar\'e inequality) in the continuum $\R^d$, and also \cite{daipra} for an approach through the entropic inequality. In all these papers, the objective is to find constants which are independent of $\Lambda$ and of the boundary condition. In a recent work \cite{posta}, Dai Pra and Posta established the entropic inequality with constant $\rho_0 \geq 1- \Vert \lambda \Vert _{\infty} \, \varepsilon (\beta)$, under the following Dobrushin-type uniqueness condition:
\begin{eqnarray}
\label{eq:epsilon}
\varepsilon (\beta) \, := \, \sum_{x\in \Z ^d} \left( 1-e^{-\beta \phi (x)} \right) & < & \frac{1}{\Vert \lambda \Vert _{\infty}} .
\end{eqnarray}
Here $\Vert \cdot \Vert_\infty$ denotes the supremum norm of a bounded function on $\Z ^d$. Note that assumption (\ref{eq:epsilon}) will be verified as soon as $\beta$ is small enough, i.e. the temperature of the system is sufficiently high. Therefore, if we choose for some $\kappa >1$ the test function $V(\eta) := \kappa ^{\sum_{x\in \Lambda} \eta_x}$ which is in $L^1(\mu)$, then an observable $f$ belongs to the class $\L_V (a,b)$ if and only if
\begin{eqnarray*}
\Gamma (f,f)(\eta) & \leq & \frac{a (\kappa -1)}{\kappa} \, \sum_{x\in \Lambda} \left( \eta_x - \kappa \, \lambda (x) \, e^{-\beta D_x ^+ H(\eta)} \right) + b , \quad \eta \in \N ^\Lambda ,
\end{eqnarray*}
as in the context of birth-death processes above. Thus the Gaussian-exponential concentration estimate of theorem~\ref{theo:result} applies under these observables. Finally, we have $D_x ^+ H(\eta) \geq 0$ because $\phi$ is non-negative and if $f\in \Lip (\N ^\Lambda)$ then one may choose
$$
a := \frac{\kappa}{2(\kappa-1)} \quad \mbox{and} \quad b := \frac{(1+\kappa)}{2} \, \sum_{x\in \Lambda} \lambda (x).
$$
In particular if $\lambda$ is assumed to be summable on $\Z^d$, then $b$ no longer depends on the box $\Lambda$. \vspace{0.5cm} \\
\textbf{Acknowledgments.} The authors are grateful to Paolo Dai Pra and G. Posta for useful discussions about the statistical mechanics part. They thank also the ANR Project EVOL for financial support.

\end{document}